\documentclass{elsarticle}

\usepackage{sets,defs,words,letters,url}

\newtheorem{theorem}{Theorem}
\newtheorem*{theorem*}{Theorem}

\newtheorem*{corollary*}{Corollary}
\newtheorem*{proposition}{Proposition}
\newtheorem*{lemma}{Lemma}

\theoremstyle{definition}
\newtheorem*{definition}{Definition}

\theoremstyle{remark}
\newtheorem*{remark}{Remark}

\begin{document}

\title{The deformation theory of sheaves of \\ commutative rings}
\date{\today}
\author{Jonathan Wise\footnote{Present address:  
Department of Mathematics,
Stanford University,
450 Serra Mall,
Building 380,
Stanford, CA, USA 94305 \\ \indent \hphantom{$^1$} E-mail address:  \url{jonathan@math.stanford.edu}}}
\address{Department of Mathematics,
The University of British Columbia,
Room 121, 1984 Mathematics Road,
Vancouver, B.C., Canada V6T 1Z2}

\begin{abstract}
  We define a sheaf of abelian groups whose cohomology is represented by the cotangent complex, permitting a rapid introduction to the theory of the cotangent complex in the same generality as it was defined by Illusie, but avoiding simplicial methods.  We show how obstructions to some standard deformation problems arise as the classes of torsors under and gerbes banded by this sheaf.  This generalizes results of Quillen, Rim, and Gaitsgory.
\end{abstract}

\maketitle

\section{Introduction}
\label{sec:intro}

If $f : X \rightarrow Y$ is a \emph{smooth} morphism of schemes, and $Y'$ is a square-zero extension of $Y$ such that the ideal of $Y$ in $Y'$ is $\cO_Y$, then the square-zero extensions $X'$ of $X$ over $Y'$ are obstructed by a class $\omega \in H^2(X, T_{X/Y})$, where $T_{X/Y}$ is the relative tangent bundle.  Should this class vanish, deformations form a torsor under $H^1(X, T_{X/Y})$ and automorphisms of any fixed deformation are in bijection with $H^0(X, T_{X/Y})$.

This may be explained succinctly by allowing $X$ to vary in the definition of the deformation problem:  the deformations of open subsets of $X$ form a gerbe over $X$, banded by $T_{X/Y}$.  The obstruction $\omega$ is then the class of this gerbe.  If $\omega$ vanishes, the gerbe is trivial, meaning it is isomorphic to the classifying stack of $T_{X/Y}$-torsors, and then the statements regarding isomorphism classes and automorphisms come from the cohomological classification of torsors.  

If we relax the hypothesis that $f$ be smooth, the above argument fails, but a similar description of obstructions, deformations and automorphisms persists \cite[Th\'eor\`eme~2.1.7]{Illusie}:  letting $\bL_{X/Y}$ denote the cotangent complex of $X$ over $Y$, there is an obstruction $\omega \in \Ext^2(\bL_{X/Y}, \cO_X)$ whose vanishing is equivalent to the existence of a deformation; if a deformation exists, all deformations form a torsor under $\Ext^1(\bL_{X/Y}, \cO_X)$ and automorphisms of any given solution are in bijection with $\Ext^0(\bL_{X/Y}, \cO_X)$.

The strong resemblance between this result and the one obtained in the smooth case hints that there may still be a relationship between deformations and banded gerbes.  Moreover, it is only the ``local triviality'' aspect of a gerbe that fails to apply in the non-smooth case:  the solutions to the deformation problem still form a pseudo-gerbe banded by $T_{X/Y}$, in the sense that isomorphisms between any two solutions form a pseudo-torsor under $T_{X/Y}$, but solutions are not guaranteed to exist locally, and pairs of solutions are not guaranteed to be locally isomorphic.  In other words, the failure of the ``gerbe argument'' in the non-smooth case may be attributed to the fact that the deformation problem is not locally trivial in the Zariski topology on $X$.

All of this suggests that Illusie's result may be interpreted in terms of gerbes if we can find a topology finer than the Zariski topology in which the deformation problem becomes locally trivial.  For affine schemes, such a topology was defined, apparently simultaneously, by Quillen~\cite{Q} and Rim~\cite[VI.3]{sga7-1}.\footnote{Quillen attributed his definition to Grothendieck with a pre-publication reference to SGA4 that I could not trace.  It may be that Quillen was only crediting Grothendieck with the idea of using the topology generated by universal effective epimorphisms, and not specifically for the definition of the cotangent cohomology in this way.  On the other hand, Quillen's topology is identical to Rim's, so it may also be that Quillen intended to refer to Rim's expos\'e in SGA7.}  Rim speculated \cite[VI.3.16]{sga7-1} that it might be possible to define an analogous topology for arbitrary schemes, and Quillen apparently made such a definition~\cite{Ill-rec} but never published it.  In~\cite{Ga}, Gaitsgory defined a topology on the category of associative algebras on a scheme and showed that it is fine enough to find local trivializations of deformation problems associated to quasi-coherent algebras.  As Gaitsgory notes \cite[Section~0.4]{Ga}, his methods may be adapted easily to the commutative case, where they can be used to treat the relative deformation theory of a scheme that is affine over the base.  We note, however, that if $f : X \rightarrow Y$ is a morphism of schemes then $\cO_X$ is, in general, not quasi-coherent as a $f^{-1} \cO_Y$-algebra, so Gaitsgory's results do not apply directly to the deformation theory of schemes.

The introduction of banded gerbes to explain the obstruction to the existence of algebra extensions is due to Gaitsgory~\cite{Ga1}.  That deformations, when they exist, can be viewed as torsors was observed by Quillen \cite[Proposition~2.4~(iv)]{Q}.

In this paper, we will define a new topology on the category of all commutative rings in a topos and show that it is fine enough to trivialize the standard deformation problems about commutative rings, but is still coarse enough to glue their solutions.  We obtain cohomological obstructions to the existence of solutions to these problems and a cohomological description of the solutions, should they exist.  As we explain in Section~\ref{sec:summary}, this can be used to apply the ideas of Gaitsgory, Quillen and Rim to the deformation theory of schemes.  We also compare our approach to cotangent cohomology with Illusie's, showing that our obstruction groups agree with his.  We will complete this comparison in \cite{coh-rings} by showing that the obstruction \emph{classes} agree as well.

The ideas in this paper may be applied easily to similar deformation problems of other algebraic objects.  We leave these applications to the reader for now.  We hope to explain some of them (such as stable maps and sheaves of modules) in future work.  In \cite{coh-rings}, we will explain how the theory developed here can be used in place of sheaves of simplicial commutative rings to develop the standard properties of the cotangent complex \cite[II.2]{Illusie}.

\section{Summary of results}
\label{sec:summary}

In this section we shall state our results in the context of schemes in order to give them a geometric appearance; statements in the generality of ringed topoi follow in the body of the text.  In order to deduce the statements about schemes given in this section from the algebraic statements that follow, one need only observe that infinitesimal extensions of a scheme are equivalent to infinitesimal extensions of their structure sheaves, in either the Zariski or \'etale topology.

Suppose that $f : X \rightarrow Y$ is a morphism of $Z$-schemes, and consider the problem of extending $f$ to a fixed square-zero extension $X'$ of $X$ over $Z$ with ideal $J$ (cf.\ \cite[Probl\`eme~III.2.2.1.2]{Illusie}).  Diagramattically, we are attempting to complete the commutative diagram of solid arrows
\begin{equation*}
  \xymatrix{
    X \ar[r] \ar[d] & Y \ar[d] \\
    X' \ar[r] \ar@{-->}[ur] & Z
  }
\end{equation*}
with a dashed arrow making both triangles commute.  In Section~\ref{sec:topology}, we shall define a site $g^{-1} \cO_Z\uAlg/ f^{-1} \cO_Y$, which we show in Section~\ref{sec:def-hom} is fine enough to ensure that the deformation problem is locally trivial.  In the statement of the theorem, we abbreviate the name of the site to $\cO_Z\uAlg/\cO_Y$.
\begin{theorem*}
  The extensions $\Hom^X(X', Y)$ of $f$ to $f' : X' \rightarrow Y$ form a torsor on $\cO_Z\uAlg/\cO_Y$ under the sheaf of abelian groups $\uDer_{\cO_Z}(\cO_Y, J)$ (defined in Section~\ref{sec:def-hom}).  The class of this torsor in $H^1(\cO_Z\uAlg/\cO_Y, \uDer_{\cO_Z}(\cO_y,J))$ obstructs the existence of a lift.  Provided that the obstruction vanishes, all lifts form a torsor under $\Der_{\cO_Z}(\cO_Y,J)$.
\end{theorem*}

Now consider a morphism of schemes $f : X \rightarrow Y$ and a fixed extension $Y'$ of $Y$ with ideal~$I$.  Also assume given a homomorphism $\varphi : f^\ast I \rightarrow J$ for some quasi-coherent sheaf $J$ on $X$.  We search for a completion of the diagram
\begin{equation} \label{eqn:2}
  \xymatrix{
    X \ar@{-->}[r] \ar[d] & X' \ar@{-->}[d]^{f'} \\
    Y \ar[r] & Y'
  }
\end{equation}
in which $X'$ is a square-zero extension of $X$ by the ideal $J$ and the induced morphism ${f'}^\ast I \rightarrow J$ coincides with $\varphi$ (cf.\ \cite[Probl\`eme~III.2.1.2.1]{Illusie}).  We show in Section~\ref{sec:def-theory} that this deformation problem also becomes locally trivial in $f^{-1} \cO_Y\uAlg/\cO_X$.  Once again abbreviating the name of the site to $\cO_Y/\uAlg/\cO_X$, we obtain
\begin{theorem*}
  The completions of Diagram~\eqref{eqn:2} form a gerbe on $\cO_Y\uAlg/\cO_X$ banded by $\uDer_{\cO_Y}(\cO_X, J)$.  The class in $H^2\bigl(\cO_Y\uAlg/\cO_X, \uDer_{\cO_Y}(\cO_X, J)\bigr)$ of this gerbe obstructs the existence of a solution to this problem.  If the obstruction vanishes, solutions form a torsor under $H^1\bigl(\cO_Y\uAlg/\cO_X, \uDer_{\cO_Y}(\cO_X,J)\bigr)$, and the automorphisms group of any solution is $H^0\bigl(\cO_Y\uAlg/\cO_X, \uDer_{\cO_Y}(\cO_X,J)\bigr)$.
\end{theorem*}

In the case where $I = 0$ and $Y = Y'$, an extension always exists, so this implies
\begin{corollary*}
  There is an equivalence of categories between the category of completions of the diagram
  \begin{equation*} \xymatrix{
      X \ar@{-->}[r] \ar[d] & X' \ar@{-->}[dl] \\
      Y,
    }
  \end{equation*}
  by a scheme $X'$ that is a square-zero extension of $X$ with ideal $J$ and the category of torsors on $\cO_Y\uAlg/\cO_X$ under the sheaf of abelian groups $\uDer_{\cO_Y}(\cO_X, J)$.  Isomorphism classes are in bijection with $H^1\bigl(\cO_Y\uAlg/\cO_X, \uDer_{\cO_Y}(\cO_X, J)\bigr)$ and automorphisms of any object are in bijection with $H^0\bigl(\cO_Y\uAlg/\cO_X, \uDer_{\cO_Y}(\cO_X, J)\bigr)$.
\end{corollary*}

How do these obstructions compare to those defined by Illusie?  In Section~\ref{sec:free}, we prove
\begin{theorem*}
  The cotangent complex $\bL_{X/Y}$ \cite[II.1.2.3]{Illusie} represents cohomology of the sheaves $\uDer_{\cO_Y}(\cO_X, J)$ on the site $\cO_Y\uAlg/\cO_X$, in the sense that
  \begin{equation*}
    \Ext^p(\bL_{X/Y}, J) = H^p\bigl(\cO_Y\uAlg/\cO_X, \uDer_{\cO_Y}(\cO_X, J)\bigr) .
  \end{equation*}
\end{theorem*}
\noindent This shows that our obstruction groups are the same as Illusie's.  In \cite{coh-rings}, we show that the obstruction \emph{classes} agree with Illusie's by the identification above.

\section{Review of topologies, sites, and topoi}
\label{sec:topos-review}

If $X$ is an object of a category $C$, a sieve of $X$ is a subfunctor of the functor of points of $X$.  It is frequently preferable to view a sieve as a fibered subcategory over $C$ of the category $C/X$ of objects of $C$ over $X$.  We shall pass back and forth between these perspectives without comment.  See \cite[I.4]{sga4-1} for more above sieves (French: \emph{cribles}).

The sieve of $X$ generated from a collection $S$ of maps $Y \rightarrow X$ is (as a fibered category) the collection of all $Z \rightarrow X$ that factor through some $Y \rightarrow X$ in $S$.

A topology on a category $C$ is a collection $J(X)$ of sieves of $X$ for each object $X$ of $C$.  These are generally called the \emph{covering sieves} of the topology.  The collections are required to comply with the conditions T1 (change of base), T2 (local character), and T3 (inclusion of the sieves generated by identity maps) of \cite[D\'efinition~II.1.1]{sga4-1}.  We shall usually describe the sieves in each $J(X)$ by giving generators.  These generators are called \emph{covering families}.

Every category possesses a \emph{canonical topology} \cite[D\'efinition~II.2.5]{sga4-1}, the finest in which all representable functors are sheaves.

If $C$ is a category with a topology, a family $S$ of objects of $C$ is called a collection of \emph{topological generators} of $C$ if every object $X$ of $C$ is covered by a sieve generated by maps from objects of $S$ to $X$ \cite[D\'efinition~II.3.0.1]{sga4-1}.  We shall depart from the definitions \cite[I.1.1.5]{sga4-1} and \cite[III.2.1]{MM} and call $C$, together with its topology, a \emph{site} if $C$ possesses a set of topological generators.

If $C$ is a category with a topology, a functor $F : C^\circ \rightarrow \Sets$ is called a sheaf if the natural map
\begin{equation*}
  F(X) \rightarrow \Hom(R, F)
\end{equation*}
is a bijection for every $X$ in $C$ and for every covering sieve $R$ of $X$ \cite[D\'efinition~II.2.1]{sga4-1}.  Here, $\Hom(R, F)$ refers to maps of presheaves.  If $C_0$ is the full subcategory of $C$ spanned by a collection of generators of $C$ with its induced topology, then the restriction map identifies the category of sheaves on $C$ with the category of sheaves on $C_0$.  If $C$ is a site then $C$ is generated topologically by a small subcategory, which ensures that the sheaves on $C$ form a category, which is frequently denoted $\tC$.

Any presheaf on a site $C$ has an associated sheaf.  If $X$ is an object of $C$ then the associated sheaf of the functor represented by $X$ is denoted $\epsilon(X)$.  This determines a functor $\epsilon : C \rightarrow \tC$.

A site $C$ is called a \emph{topos} if its topology is the canonical one and the functor $\epsilon : C \rightarrow \tC$ described above is an equivalence.  Once again, our definition is slightly different from \cite[D\'efinition~IV.1.1]{sga4-1}; it is equivalent to \cite[Definition~III.4.3]{MM}.  Every topos $E$ has a final object, which we will denote by the same letter $E$.

A morphism of topoi $f : C \rightarrow C'$ is a pair of functors $f^\ast : C' \rightarrow C$ and $f_\ast : C \rightarrow C'$ such that $f_\ast$ is right adjoint to $f^\ast$ and $f^\ast$ is exact.  Every topos admits an essentially unique morphism to the topos $\Sets$.

We shall depart again from \cite[IV.4.9.1]{sga4-1} in the definition of a morphism of sites, and declare that a morphism between sites is a morphism between their associated topoi.  This extrinsic definition can be made intrinsic (see e.g.\ \cite[Theorem~VII.10.1]{MM}), but we will be content to recall here that a morphism of sites can be induced from a \emph{cocontinuous functor} $C \rightarrow C'$ or from a left exact \emph{continuous functor} $C' \rightarrow C$.  

A functor $f : C \rightarrow C'$ between sites is cocontinuous if, for any $X \in C$ and any covering sieve $R$ of $f(X)$, the collection of all $Y \rightarrow X$ such that $f(Y) \in R$ is a covering sieve \cite[D\'efinition~III.2.1, Proposition~III.2.3]{sga4-1}.  A left exact functor $g : C' \rightarrow C$ is continuous if and only if it takes covering families to covering families \cite[D\'efinition~III.1.1, Proposition~III.1.3, Proposition~III.1.6]{sga4-1}.  If $f$ is left adjoint to $g$ then cocontinuity of $f$ coincides with continuity of $g$ \cite[Proposition~III.2.5]{sga4-1}.

By \cite[Th\'eor\`eme~1.10.1]{T}, every topos has enough injective sheaves of abelian groups.  This permits us to define the derived functors of $f_\ast$ (applied to sheaves of abelian groups) for any morphism of topoi $f : C \rightarrow C'$.  In the case where $C' = \Sets$, these derived functors are denoted $F \mapsto H^p(C, F)$.

Suppose that $A \xrightarrow{u} B \xrightarrow{v} C$ is a sequence of morphisms of topoi.  Since $u_\ast$ has an exact left adjoint it preserves injectives and we obtain a spectral sequence for the composition of functors \cite[Th\'eor\`eme~2.4.1]{T}
\begin{equation*}
  R^p v_\ast \big( R^q u_\ast F \big) \Rightarrow R^{p+q}(vu)_\ast F 
\end{equation*}
for any sheaf of abelian groups on $A$.  In terms of derived categories, we have $R v_\ast \circ R u_\ast = R (vu)_\ast$.  In the case where $C = \Sets$ the spectral sequence above specializes to
\begin{equation*}
  H^p(B, R^q u_\ast F) \Rightarrow H^{p+q}(A, F) .
\end{equation*}
See \cite[Expos\'e~V]{sga4-2} for more about the Cartan--Leray spectral sequences.

If $X_\bullet$ is a \emph{hypercover} of $E$  \cite[V.7.3.1]{sga4-2} then there is a spectral sequence \cite[V.7.4.0.3]{sga4-2}
\begin{equation*}
  H^p (X_\bullet, \cH^q(F)) \Rightarrow H^{p+q}(E,F) .
\end{equation*}
More specifically, if $I^\bullet$ is a resolution of $F$ by sheaves that are acyclic for each of the $X_p$, then the double complex $\Gamma(X_\bullet, I^\bullet)$ computes the cohomology of $F$.  We refer the reader to \cite[V.7]{sga4-2} for more details.

\section{Review of torsors and gerbes}
\label{sec:tandg}

We restrict attention to torsors and gerbes under abelian groups.  Suppose that $E$ is a topos and $G$ is a sheaf of abelian groups on $E$.  A \emph{$G$-torsor} is a sheaf $F$ on $E$ with an action $a : G \times F \rightarrow F$ of $G$ such that
\begin{enumerate}
\item (pseudo-torsor) the map $(a, p_2) : G \times F \rightarrow F \times F$ is an isomorphism of sheaves, and
\item (local triviality) $F$ covers $E$.
\end{enumerate}
If $F$ satisfies only the first condition then $F$ is called a \emph{pseudo-torsor under $G$}.  The second condition says that, locally in $E$, the sheaf $F$ admits a section.  Since a section of a torsor trivializes it, the second condition says that $F$ is locally isomorphic to $G$ as a sheaf with $G$-action.

\begin{theorem*}[cf.\ {\cite[Remarque~III.3.5.4]{Gir}}]
  Isomorphism classes of torsors under an abelian group $G$ are in bijection with $H^1(E, G)$.  Isomorphisms between any two torsors form a pseudo-torsor under $H^0(E, G)$.  In particular, automorphisms of a torsor are in canonical bijection with $H^0(E,G)$.
\end{theorem*}

A \emph{gerbe on $E$ banded by $G$} is a stack $F$ \cite[D\'efinition~1.2.1]{Gir} on $E$ with an action of $G$ on the morphism sheaves of $G$ that is compatible with composition and satisfies
\begin{enumerate}
\item (pseudo-gerbe) for any two sections $\xi, \eta$ of $F$ over $U$, the sheaf $\uIsom(\xi, \eta)$ is a $G$-pseudo-torsor on $U$,
\item (local triviality for morphisms) any two sections of $F$ over $U$ are locallly isomorphic,
\item (local triviality for objects) $F$ covers $E$.
\end{enumerate}
If $F$ satisfies only the first condition, we call $F$ a pseudo-gerbe banded by $G$.  The second condition ensures that $\uIsom(\xi, \eta)$ forms a $G_U$-torsor for each pair of sections $\xi, \eta \in \Gamma(E/U, F)$.  As in the case of torsors, the final condition means that sections of $F$ exist locally in $E$.  Since a section of $F$ induces an equivalence of banded gerbes between $F$ with $BG$, the classifying stack of $G$-torsors, we say that a banded gerbe is locally isomorphic to $BG$.

\begin{theorem*}[cf.\ {\cite[Th\'eor\`eme~IV.3.4.2]{Gir}}]
  Equivalence classes of gerbes on $E$ banded by an abelian group $G$ are in bijection with $H^2(E, G)$.  If $\sF$ is a gerbe banded by $G$ then sections of a $\sF$ form a pseudo-torsor under $H^1(E, G)$.  Isomorphisms between any two sections of $\sF$ form a torsor under $H^0(E, G)$.  In particular, automorphisms of any section of $\sF$ are in canonical bijection with $H^0(E, G)$.
\end{theorem*}

\section{The topos of commutative rings}
\label{sec:topology}

All rings and algebras are commutative and unital.

Now let $(E,A)$ be a ringed topos \cite[IV.11.1.1]{sga4-1}, and let $B$ an $A$-algebra.  Let $A\uAlg(E)/B$ (or $A\uAlg/B$ for short) be the category of pairs $(U, C)$ where $U \in E$ and $C$ is an $A_U$-algebra with a map to $B_U$.  A morphism of $A\uAlg/B$ from $(U_1, C_1)$ to $(U_2, C_2)$ is a map $f : U_1 \rightarrow U_2$ of $E$ and a map $C_1 \rightarrow f^\ast C_2$ of $A_{U_1}$-algebras commuting with the projections to $B_{U_1}$.

\begin{definition}
  A family of maps $(U_i, C_i) \rightarrow (U, C)$ in $A\uAlg/B$ is covering if, for any $V \rightarrow U$ and any finite set of sections $\Lambda \subset \Gamma(V, C)$, there exists, locally in $V$, a map $V \rightarrow U_i$ for some $i$ and a lift of $\Lambda$ to $\Gamma(V, C_i)$.
\end{definition}

To understand this topology, it may be helpful to consider the case where $E$ is a point.  In that case, $A \rightarrow B$ is a ring homomorphism, and a family of $A$-algebra maps $C_i \rightarrow C$ over $B$ is considered to be covering if every finite set of elements of $C$ can be lifted to some $C_i$.

Now let $E$ be an arbitrary topos and suppose that $R$ is a sieve of $(U,C)$ in $A\uAlg(E)/B$.  If $V \rightarrow U$ is a map of $E$ and $\Lambda \subset \Gamma(V, C)$ is finite, let $(V, A_V[\Lambda]) \rightarrow (U, C)$ be the induced map of $A\uAlg/B$.  Let $Q(\Lambda)$ the collection of all $W \rightarrow V$ such that it is possible to complete the diagram
\begin{equation*}
  \xymatrix{
    (W, A_W[\Lambda]) \ar[r] \ar@{-->}[d] & (V, A_V[\Lambda]) \ar[d] \\
    (U', C') \ar[r] & (U, C)
  }
\end{equation*}
with $(U', C')$ in $R$.  Then for $R$ to be a covering sieve means that $Q(\Lambda)$ is a covering sieve of $V$ in $E$ for every $Q(\Lambda)$ arising as above.

It is immediate from this description of the topology that any $(U, C)$ in $A\uAlg/B$ is covered by the collection of all $(V, A_V[S]) \rightarrow (U,C)$ where $S$ is an arbitrary finite set.  Therefore the pairs $(V, A_V[S])$ generate the topology of $A\uAlg/B$ and we are free to say that in this topology, any $A$-algebra is locally a finitely generated polynomial ring.  Furthermore, we obtain a set of topological generators for $A\uAlg/B$ by taking the collection of all $(V, A_V[S])$ such that $V$ lies in a set of topological generators for $E$.

\begin{remark}
  This topology is slightly more complicated than the ones used by Gaitsgory, Quillen, and Rim.  Analogues of those topologies would work here, but this topology has a technical advantage in its possession of a set of topological generators.  This permits us to make use of topoi without making recourse to universes.
\end{remark}

\section{Deformation of homomorphisms}
\label{sec:def-hom}

Let $A$ be a sheaf of rings on $E$ and $B \rightarrow C$ a homomorphism of $A$-algebras.  Suppose that $C'$ is a square-zero extension of $C$, as an $A$-algebra, with ideal $J$.  Consider the problem of lifting the homomorphism $B \rightarrow C$ to a map $B \rightarrow C'$ (cf.\ \cite[Probl\`eme~III.2.2.1.1]{Illusie}).

Putting $B' = C' \fp_C B$, this problem immediately reduces to that of finding a section of $B'$ over $B$.  We denote the set of such sections $\Hom^A_B(B, B')$.

The difference between any two sections of $B'$ over $B$ is an $A$-derivation from $B$ into $J$.  Denoting the set of all such derivations by $\Der_A(B,J)$ we may say that $\Hom^A_B(B, B')$ is a $\Der_A(B,J)$-pseudo-torsor.  If $B$ is allowed to vary we obtain a sheaf $\uDer_A(B,J)$ on $A\uAlg/B$, represented by the object $(E, B + J)$, where $B + J$ is the trivial square-zero extension of $B$ by $J$.  Then $B'$ represents a $\uDer_A(B,J)$-pseudo-torsor on $A\uAlg/B$.

In fact, $B'$ is a $\uDer_A(B,J)$-torsor on $A\uAlg/B$, since the map $(E, B') \rightarrow (E, B)$ is covering in $A\uAlg/B$.  The cohomological classification of torsors now implies
\begin{theorem} \label{thm:1}
  Let $\omega \in H^1\big(A\uAlg/B, \: \uDer_A(B,J)\big)$ be the class of $B'$ as a torsor under $\uDer_A(B,J)$.  Then $\omega = 0$ if and only if $B'$ admits a section as an $A$-algebra over $B$.  In that case, the sections form a torsor under
  \begin{equation*}
    H^0\big(A\uAlg/B, \: \uDer_A(B, J)\big) = \Der_A(B,J) .
  \end{equation*}
\end{theorem}

\section{Extensions of algebras}
\label{sec:exal}

Suppose that $B$ is an $A$-algebra in the topos $E$ and $J$ is a $B$-module.  Let $\Exal_A(B,J)$ be the category of square-zero $A$-algebra extensions of $B$ with ideal $J$.  These categories fit together into a fibered category $\uExal_A(B,J)$ over $A\uAlg/B$.  We saw in Theorem~\ref{thm:1} that any $B' \in \Exal_A(B,J)$ represents a $\uDer_A(B,J)$-torsor on $A\uAlg/B$, so we obtain a fully faithful functor
\begin{equation} \label{eqn:3}
  \uExal_A(B,J) \rightarrow B \uDer_A(B,J)
\end{equation}
from $\uExal_A(B,J)$ to the classifying stack of $\uDer_A(B,J)$-torsors on $A\uAlg/B$.

\begin{lemma}
  Every $\uDer_A(B,J)$-torsor is representable by a square-zero $A$-algebra extension of $B$ by~$J$.
\end{lemma}

Suppose that $P$ is a $\uDer_A(B,J)$-torsor.  Let $B'$ be the sheaf on $E$ whose sections over $U$ are pairs $(b, b')$ where $b$ is a section of $B$ over $U$, corresponding to a map $b : A_U[x] \rightarrow B$ (denoted by the same letter), and $b'$ is a section of $b^\ast P$ over $(U, A_U[x])$.

We can give $B'$ a ring structure as follows.  Suppose $(b, b')$ and $(c, c')$ are two sections of $B'$.  Choose a cover $R$ of $(E, B)$ over which $P$ is trivial.  Then $E$ can be covered by objects $U$ such that there is a $(U, C) \in R$ and both $b_U$ and $c_U$ lift to $\Gamma(U, C)$.  There is therefore a map 
\begin{equation*}
  (b_U, c_U) : A_U[x,y] \rightarrow C
\end{equation*}
and $(b_U, c_U)^\ast P_C$ is trivial because $P_C$ is trivial.  A trivial $\uDer_A(B,J)$-torsor is certainly representable (by $B + J$), so let $B'_C$ be an extension of $C$ by $J$ representing $(b_U, c_U)^\ast P_C$.  We are given maps $b'_U : A_U[x] \rightarrow B'_C$ and $c'_U : A_U[y] \rightarrow B'_C$ over $C$.  Since $B'_U$ is a ring, these extend uniquely to a map $(b'_U, c'_U) : A_U[x,y] \rightarrow B'_C$ over $C$.  Restricting this, via the maps $A_U[x] \cong A_U[x+y] \rightarrow A_U[x,y]$ and $A_U[x] \cong A_U[xy] \rightarrow A_U[x,y]$, yields sections that we will denote $b'_U + c'_U$ and $b'_U c'_U$ of $(b + c)^\ast P_C$ and $(bc)^\ast P_C$, respectively.  These give us sections $(b_U + c_U, b'_U + c'_U)$ and $(b_U c_U, b'_U c'_U)$ of $B'$ over $U$.

The uniqueness of the constructions above implies that they patch together to give a ring structure on $B'$ over $E$, which makes $(E,B')$ an object of $A\uAlg/B$.  The verifications of commutativity, associativity, etc.\ are essentially the same, using a trio of sections of $C$ instead of a pair.  To check that $B'$ represents $P$, one only needs to produce an isomorphism between $P$ and the object represented by $B'$ under the assumption that $P$ admits a section and $B$ is a finitely generated polynomial algebra over $A$, since the pairs $(U, A_U[S])$ such that $P_U$ is trivial and $S$ is finite generate the topology of $A\uAlg/B$.  Under these assumptions, the construction clearly provides the isomorphism, and this varies in a functorial way with free $A$-algebras $B$.  \qed

\begin{remark}
  Let $F$ be the fibered category of pairs $(U, B \rightarrow C)$ where $B \rightarrow C$ is an $A$-algebra morphism over $U$ and morphisms are commutative squares.  The projection $F \rightarrow A\uAlg$ sending the object above to $(U, C)$ makes $F$ into a fibered category over $A\uAlg$ and the proof of the lemma demonstrates that $F$ is a stack over $A\uAlg$.
\end{remark}

The cohomological classification of torsors now gives us
\begin{theorem} \label{thm:3}
  The functor~\eqref{eqn:3} is an equivalence.  Isomorphism classes in the category $\Exal_A(B,J)$ are therefore in bijection with $H^1(A\uAlg/B, \uDer_A(B,J))$.
\end{theorem}

\section{Deformation of algebras}
\label{sec:def-theory}

Let $A$ be a sheaf of rings on $E$ and $B$ an $A$-algebra.  Suppose that $A'$ is an extension of $A$ with ideal $I$ (not necessarily square-zero).  Fix a $B$-module $J$ and an $A \rightarrow B$ homomorphism $\varphi : I \rightarrow J$.  Define $\Def_A(A', B, \varphi)$ to be the category of completions of the diagram
\begin{equation*}
  \xymatrix{
    0 \ar[r] & I \ar[r] \ar[d] & A' \ar[r] \ar@{-->}[d] & A \ar[r] \ar[d] & 0 \\
    0 \ar[r] & J \ar@{-->}[r] & B' \ar@{-->}[r] & B \ar[r] & 0 
  }
\end{equation*}
by an extension $B'$ of $B$ by $J$ (cf.\ \cite[Probl\`eme~2.1.2.1]{Illusie}).    Allowing $B$ to vary, these categories fit together into a fibered category $\uDef_A(A', B, \varphi)$ over $A\uAlg/B$.  (Note that the special case $I = 0$ recovers $\Exal_A(B,J) = \Def_A(A,B,0)$.)

If $B'$ and $B''$ are any objects in $\Def_A(A',B,\varphi)$ then the difference between any two isomorphisms between $B'$ and $B''$ is a derivation $B \rightarrow J$, i.e., an element of $\Der_A(A',B,\varphi)$.  It follows from the lemma of Section~\ref{sec:exal} that $\Der_A(A',B,\varphi)$ is a stack, so this tells us that $\Def_A(A',B,\varphi)$ is a pseudo-gerbe banded by $\Der_A(B,J)$, and $\uDef_A(A',B,\varphi)$ is a pseudo-gerbe over $A\uAlg/B$, banded by $\uDer_A(B,J)$.  In fact, we have

\begin{proposition}
  $\uDef_A(A',B,\varphi)$ is a gerbe over $A\uAlg/B$, banded by $\uDer_A(B,J)$.
\end{proposition}

We must check that $\uDef_A(A',B,\varphi)$ admits a section locally in $A\uAlg/B$, and that any two sections are locally isomorphic.  Since $B$ is locally a polynomial algebra over $A$, it's sufficient for the local existence to show that $\Def_A(A',B,\varphi)$ in that case.  But then we could take $B' = A'[S] \cp_{I[S]} J$, the extension obtained by pushing out the extension $A'[S]$ of $A[S]$ by $I[S]$ by the canonical map $I[S] \rightarrow J$.

To prove that any two sections $B'$ and $B''$ of $\uDef_A(A',B,\varphi)$ are locally isomorphic, we shall construct their difference $B''' = B'' - B'$ and show that it is a trivial extension of $B$ as an algebra over $A' - A' = A + I$.  Before making this precise, note that we may replace $A'$ by $A'/I^2$, and therefore assume that $I^2 = 0$, without changing the deformation problem.  The ring $B'' \fp_B B'$ is an extension of $B$ by the ideal $J \times J$ and there is a morphism of exact sequences,
\begin{equation*}
  \xymatrix{
    0 \ar[r]  &I \times I \ar[r] \ar[d]  &A' \fp_A A' \ar[r] \ar[d]  &A \ar[r] \ar[d] &0 \\
    0 \ar[r] &J \times J \ar[r] & B' \fp_B B'' \ar[r] & B \ar[r] & 0 .
  }
\end{equation*}
We push out these sequences via the difference maps $I \times I \rightarrow I$ and $J \times J \rightarrow J$ sending $(x,y)$ to $x - y$.  This yields a map of exact sequences
\begin{equation*}
  \xymatrix{
    0 \ar[r] & I \ar[r] \ar[d] & A + I \ar[r] \ar[d] & A \ar[r] \ar[d] & 0 \\
    0 \ar[r] & J \ar[r] & B''' \ar[r] & B \ar[r] & 0 .
  }
\end{equation*}
Note that $B''$ can be recovered functorially from $B'$ and $B'''$ by an addition procedure inverse to the difference procedure just executed.  Thus to show that $B'$ and $B''$ are locally isomorphic, as extensions of $B$ compatible with $A'$, it is equivalent to show that $B'''$ is locally isomorphic to the trivial extension $B + J$ of $B$, with its trivial $(A + I)$-algebra structure.

The $(A+I)$-algebra structure of $B'''$ is determined by the $A$-algebra structure induced from the section $A \rightarrow A + I$.  It is therefore equivalent to show that any extension of $B$ by $J$ as an $A$-algebra is locally isomorphic to $B + J$.  This was the content of Theorem~\ref{thm:1}.  \qed

The cohomological classification of banded gerbes now provides
\begin{theorem} \label{thm:2}
  Let $\omega$ be the class in $H^2\big(A\uAlg/B,\: \uDer_A(B,J)\big)$ corresponding to the banded gerbe $\uDef_A(A',B,\varphi)$.  Then $\omega = 0$ if and only if $\Def_A(A',B,\varphi)$ is non-empty.
  
  In that case, $\uDef_A(A',B,\varphi)$ is isomorphic as a banded gerbe to the classifying stack of $\uDer_A(B,J)$-torsors.  Hence, isomorphism classes in $\Def_A(A',B,\varphi)$ form a torsor under the group $H^1\big(A\uAlg/B, \:\uDer_A(B,J)\big)$, and the automorphism group of any fixed object of $\Def_A(B,J)$ is canonically isomorphic to the group $\Der_A(B,J) = H^0\big(A\uAlg/B,\: \uDer_A(B,J)\big)$.
\end{theorem}

\section{The cohomology of free algebras}
\label{sec:free}

Suppose that $E$ is a topos, $A$ a sheaf of algebras on $E$, and $S$ a sheaf of sets on $E$.  Let $J$ be a sheaf of $A[S]$-modules.  We wish to compare the cohomology groups of $\uDer_A(A[S], J)$ on $A\uAlg/A[S]$ and $J_S$ on $E/S$.  We construct several sites to mediate between the $A\uAlg/A[S]$ and $E/S$.

Let $\uSets(E)$ and $\uSets^\ast (E)$, or $\uSets$ and $\uSets^\ast$ for short, be the sites whose common underlying category is the category of pairs $(U, T)$ where $U$ is an object of $E$ and $T$ is a sheaf of sets on $U$.  A map $(U_1, T_1) \rightarrow (U_2, T_2)$ is a map $f : U_1 \rightarrow U_2$ and a map $T_1 \rightarrow f^\ast T_2$.  (From another point of view, this category is the category of arrows in $E$.)

We shall say that a family of morphisms $(U_i, T_i) \rightarrow (U,T)$ is covering in $\uSets^\ast$ if, for any $f : V \rightarrow U$ and any finite subset $\Lambda \subset \Gamma(V, f^\ast T)$, there is, locally in $V$, a factorization of $f$ through $g : V \rightarrow U_i$, for some $i$, and a lift of $\Lambda$ to $\Gamma(V, g^\ast T_i)$.  The topology on $\uSets$ is defined in the same way, except $\Lambda$ is restricted to be a $1$-element set.

\begin{remark}
  The topologies on $\uSets$ and $\uSets^\ast$ are genuinely different.  In the case where $E$ is the punctual topos (i.e., the category of sets), the category of sheaves on $\uSets(E)$ may be identified with the category of sets; the category of sheaves on $\uSets^\ast(E)$ may be identified with the category of presheaves on the category of finite sets.
\end{remark}

The topology on $\uSets$ is finer than that on $\uSets^\ast$, so there is a morphism of sites $\uSets \rightarrow \uSets^\ast$ (the identity functor is cocontinuous).  This induces a map $\Phi : \uSets / (E,S) \rightarrow \uSets^\ast / (E,S)$, for any sheaf of sets $S$ on $E$.

There is also a functor $\uSets^\ast / (E,S) \rightarrow A\uAlg/A[S]$ which sends $(U,T)$ to $(U, A_U[T])$.  This functor is left exact, and by definition, it takes covers to covers so it is continuous by \cite[Proposition~III.1.6]{sga4-1} and we get a morphism of sites $\Psi : A\uAlg/A[S] \rightarrow \uSets^\ast / (E,S)$.

Finally, we have a cocontinuous functor $\uSets(E) \rightarrow E$ sending $(U,T)$ to $i_! T$, where $i$ is the canonical morphism of topoi from $E/U$ to $E$  \cite[IV.5.1--2]{sga4-1}.  This induces a morphism of sites $\Xi : \uSets / (E,S) \rightarrow E / S$.  Putting all of these morphisms together, we obtain the following diagram of morphisms of sites.
\begin{equation*}
  \xymatrix{
    & \uSets / (E,S) \ar[dr]^{\Phi} \ar[dl]_{\Xi} & & A\uAlg/A[S] \ar[dl]_{\Psi} \\
    E / S & & \uSets^\ast / (E,S)
  } 
\end{equation*}

\begin{proposition}
   The morphisms $\Phi$ and $\Psi$  are acyclic.
\end{proposition}

First consider $\Psi : A\uAlg/A[S] \rightarrow \uSets^\ast / (E,S)$.  To see that $R^p \Psi_\ast F = 0$ for $p > 0$ is a local problem in $\uSets^\ast / (E,S)$.  We can therefore reduce to the case where $S$ is the constant sheaf associated to some finite set $S_0$, since pairs $(V,T)$ where $T$ is constant and finite generate the topology of $\uSets^\ast / E$.  We must show that if $\alpha \in H^p(A\uAlg/A[S], \Psi^\ast (U,T))$ then $\alpha$ can be trivialized on some cover of $(U,T)$.  But $\Psi^\ast (U, T) = (U, A_U[T])$ and all covering sieves of $(U, A_U[T])$ in $A\uAlg/A[S]$ are pulled back from covering sieves of the final object of $E$, hence are pulled back from covering sieves of $(E,S)$ in $\uSets^\ast / (E,S)$.  Since $\alpha$ can certainly be trivialized on some covering sieve of $(U, A_U[T])$, any covering sieve of $\uSets^\ast / (E,S)$ that pulls back via $\Psi$ to this one will trivialize $\alpha$ in $\uSets^\ast / (E,S)$.

Now consider $\Phi : \uSets / (E,S) \rightarrow \uSets^\ast / (E,S)$.  If $R$ is a covering sieve of some object $(U,T)$ of $\uSets / (E,S)$, then let $R'$ be the collection of all finite disjoint unions of objects of $R$.  If $R$ is a covering sieve of $\uSets / (E,S)$ which trivializes a cohomology class, then $R'$ is a covering sieve of $\uSets^\ast / (E,S)$ which trivializes the same cohomology class.  Therefore this map is acyclic as well.  \qed

\begin{lemma}
  There is a canonical isomorphism $\Phi_\ast (E,J) \simeq \Psi_\ast \uDer_A(A[S], J)$ for any sheaf of $A[S]$-modules $J$.
\end{lemma}

This is a matter of unwinding the definitions.  For $(U,T) \in \uSets^\ast / (E,S)$, we have
\begin{equation*}
  \Gamma\big((U,T),\: \Phi_\ast (E,J)\big) = \Gamma\big((U,T),\: (E,J)\big) = \Hom_U(T, J_U) .
\end{equation*}
On the other hand,
\begin{multline*}
  \Gamma\big((U,T), \:\Psi_\ast \uDer_A(A[S], J)\big) = \Gamma\big((U,A[T]), \:\uDer_A(A[S], J)\big)  \\
  = \Der_{A_U}(A_U[T], \:J_U) = \Hom_U(T, J_U)
\end{multline*}
using the universal property of $A[T]$.  \qed

Since $\Phi$ and $\Psi$ are acyclic, this proves that $R \Phi_\ast (E,J) = R \Psi_\ast \uDer_A(A[S], J)$.  We can therefore compute the cohomology of $\uDer_A(A[S], J)$ by computing the cohomology of $(E,J)$ on $\uSets / (E,S)$.

\begin{lemma}
  The natural map $J \rightarrow R \Xi_\ast \Xi^\ast J$ is an isomorphism.
\end{lemma}

The projection $\uSets/(E,S) \rightarrow E/S$ sending $(U, T)$ to $T$ is left exact and induces an exact left adjoint $\Xi_!$ to $\Xi^\ast$ \cite[Proposition~I.5.4~4]{sga4-1}.  Therefore, $\Xi^\ast$ preserves injectives.  Since $\Xi^\ast$ is also exact, this implies $R \Xi_\ast \Xi^\ast J = \Xi_\ast \Xi^\ast J$.  The natural map $J \rightarrow \Xi_\ast \Xi^\ast J$ is certainly an isomorphism, since $\Xi^\ast J = (E,J)$ and $\Hom_{(E,S)}\bigl((E,T), (E,J)\bigr) = \Hom_S(T,J)$.  \qed

Let $u : E/S \rightarrow E$ and $v : A\uAlg/A[S] \rightarrow E$ denote the projections.  Then the lemma implies that $R v_\ast \uDer_A(A[S],J) = R u_\ast J$.

\begin{lemma}
  There is a canonical isomorphism $R \uHom(\bZ^S, J) \simeq R u_\ast u^\ast J$ for any sheaf of abelian groups $J$ on $E$.
\end{lemma}

Since $u^\ast$ has an exact left adjoint on sheaves of abelian groups~\cite[Proposition~IV.11.3.1]{sga4-1}, both $R \uHom(\bZ^S, J)$ and $R u_\ast u^\ast J$ can be computed by taking an injective resolution of $J$ in $E$.  It's therefore sufficient to remark that $\uHom(\bZ^S, J) = \uHom(S, J)$, by definition.  \qed

Putting all of the lemmas together, we find that
\begin{equation*}
  R v_\ast \uDer_A(A[S], J) = R \uHom_{\bZ}(\bZ^S, J) = R \uHom_{A[S]}(A[S]^S, J) = R \uHom(\Omega_{A[S] / A}, J) .
\end{equation*}

\begin{corollary*}
  If $J$ is injective and $B$ is a free $A$-algebra then for every $p > 0$, the group $H^p(A\uAlg/B, \uDer_A(B,J))$ vanishes.
\end{corollary*}

Since $\Omega_{B/A}$ is functorial in $B$, this permits us to compute the cohomology of $\uDer_A(B, J)$ for any sheaf $B$ of $A$-algebras and any $B$-module $J$ using hyper-\v{C}ech cohomology.  Let $B_\bullet$ be the standard simplicial resolution of $B$ by free $A$-algebras \cite[I.1.5.5~b, II.1.2.1.1]{Illusie}.  By \cite[I.1.5.3]{Illusie}, the resolution $B_\bullet \rightarrow B$ is a homotopy equivalence on the underlying sheaves of sets.  This implies in particular that it is a hypercover, so corollary above implies that $R v_\ast \uDer_A(B,J)$ is computed by the complex $C$ with
\begin{equation*}
  C^{p,q} = v_\ast \uDer_A(B_p, J^q)
\end{equation*}
for any injective resolution $J^\bullet$ of $J$.  Taking $B_\bullet = P_A(B)$ we obtain the cotangent complex as $\bL_{B/A} = \Omega_{B_\bullet / A} \tensor_{B_\bullet} B$ by \cite[II.1.2.3]{Illusie}, and then
\begin{multline*}
\displaystyle  R v_\ast \uDer_A(B, J) \simeq v_\ast \uDer_A(B_\bullet, J^\bullet) = \uHom_{B_\bullet}(\Omega_{B_\bullet / A}, J^\bullet) \\ = \uHom_B(\Omega_{B_\bullet / A} \tensor_{B_\bullet} B, J^\bullet) \simeq R \uHom(\bL_{B/A}, J) .
\end{multline*}

This proves that the cotangent complex represents the functor sending $J$ to the cohomology of $\uDer_A(B, J)$:
\begin{theorem}
  Let $v$ denote the projection $A\uAlg(E)/B \rightarrow E$.  If $J$ is any $B$-module, then there is an isomorphism in the derived category of sheaves of $B$-modules
  \begin{equation*}
    R v_\ast \uDer_A(B, J) \simeq R \uHom(\bL_{B/A}, J) .
  \end{equation*}
  In particular, 
  \begin{equation*}
    \Ext^p(\bL_{B/A}, J) = H^p\big(A\uAlg(E)/B, \:\uDer_A(B,J)\big)
  \end{equation*}
  for all $p$.
\end{theorem}

\section{Acknowledgements}

I am happy to thank Dan Abramovich, Patrick Brosnan, Barbara Fantechi, Martin Olsson, Ravi Vakil, and Angelo Vistoli for helpful conversations about these topics.  Although a goal of this work has been to circumvent some of the technical aspects of~\cite{Illusie}, Illusie's treatise has been a continual inspiration.  Vistoli's remark that ``with uncanny regularity [deformations are] a cohomology group of a certain algebraic object, and [obstructions are] the cohomology group of the same object in one degree higher'' \cite[p.~2]{Vis} was a particular impetus for this project.  Dan Abramovich provided numerous suggestions that improved the exposition of this paper considerably (and helped me to avoid one serious error).  I am also grateful to Bhargav Bhatt for informing me of Gaitsgory's work.

This research was partly supported by NSF-MSPRF 0802951.

\bibliographystyle{amsalpha}
\bibliography{def-rings}

\end{document}